УДК 519.21

# Особенности динамики частицы когда случайные возмущения ортогональны к ее скорости

В.А. Дубко

Kyiv National University of Technologies and Design
Education and Scientific Institute of Modern Learning Technologies
Department of higher mathematics.
(01011 Kyiv, 2 Nemyrovycha-Danchenka St., building 4)
doobko2017@ukr.net

**Аннотация.** Мы исследуем свойства решений уравнения Ланжевена, когда стохастические воздействия ортогональны к скорости частицы. Процесс Винера может принимать неограниченные значения. Но для этих уравнений, существуют замкнутые притягивающие поверхности. Для таких стохастических уравнений, мы построили уравнения для плотности вероятности положения частицы в пространстве координат, которые зависят от начального вектора скорости частицы. Опираясь на нашу более раннюю работу и физический смысл коэффициентов, мы строим уравнение для диффузионного приближения исходного уравнения, и находим его решение. Показано, что, при определенных соотношениях между коэффициентов в начальном стохастическом уравнении, когда силы влияния стремятся к нулю, распределение положения частицы может быть аппроксимировано решением волнового уравнения с постоянной скоростью.

# The specificity of the particle dynamics if random perturbations are orthogonal to its velocity

**V.A. Doobko**

**Annotation.** We explore properties the solution of Langevin equation when stochastic influence is orthogonal to velocity of a particle. Wiener's process can accept unlimited values. But for these equations, the attraction surfaces exist. For these stochastic equations we have constructed the equations for density of probability in co-ordinates space that are depending an on initial vector of the particle velocity. Then, using our earlier work and physical sense of coefficients, we build the equation for diffusion approximating of the original equation, and find its solution. It is shown that when at the certain concordance between of coefficients in initial stochastic equation, and influences aspire to the zero; the position distribution of the particle can be approximated by the solution of the wave equation with constant speed.

## Введение.

В уравнении Ланжевена [1], моделирующего изменение скорости броуновской частицы, входят два слагаемые. Относительно коэффициентов этого уравнения, можно сделать определенные заключения, исходя из физических соображений.

Первое слагаемое – это торможение, вызванное осредненным действием, разницы между импульсами молекул, «набегающих» и с «догоняющих» частицу. Второе - моделирует стохастичность отклонений от осредненного влияния на скорость частицы. Для больших скоростей частицы, интенсивность отклонений от усредненного воздействия по направлению движения, уменьшается. В тоже время, при движении в пространстве размерности выше двух, интенсивность воздействий перпендикулярных к направлению движения частицы не зависит от ее скорости [2, с.46-47] (аналогия - «перемещение под дождем»). Это мотивирует рассмотрение моделей с разделением случайных воздействий на скорость движения броуновской частицы: по направлению движения и перпендикулярных ему.

Рассмотрим модель, в которой присутствуют только перпендикулярные к скорости частицы случайные воздействия:

$$dv(t) = -a(t)v(t)dt + F(t)dt + b(t)P(v(t))d\xi(t) \quad (1)$$

$$dx(t) = v(t)dt \quad (2)$$

где $a(t)$ и $b(t)$ – скалярные функции, $F(t)$ - вектор внешних детерминированных воздействий, $x(t), v(t)$ - текущее положение частицы и ее вектор скорости, $\xi(t)$ векторный случайный процесс, $P(v_t)$ - оператор ортогонализации:

$$(v(t), P(v_t)d\xi(t)) \equiv 0, \quad \forall t \geq 0 -$$

где $(,)$ - обозначение скалярного произведения.

Например, для двумерных процессов, когда $\xi(t) = \tilde{\xi}(t)(e_1, e_2)$ где $\tilde{\xi}(t)$ скаляр, $e_1, e_2$ - ортогональные единичные векторы в выбранной системе координат, то

$$P(v(t))d\tilde{\xi}(t) = \frac{1}{|v(t)|}\det\begin{pmatrix} e_1 & e_2 \\ v_1(t) & v_2(t) \end{pmatrix} d\tilde{\xi}(t),$$

или, когда у $\xi(t)$ - разные компоненты, то

$$P(v(t))d\xi(t) = \frac{v(t)}{|v(t)|^2}(v(t), d\xi(t)) - d\xi(t).$$

Деление на $|v(t)|$, $|v(t)|^2$ обеспечивает зависимость воздействия только от направления $v(t)$.

Класс таких уравнений будем называть моделями с ортогональными воздействиями. Модели (1), обладают рядом особенностей, отличающих их как от моделей со случайным, скачкообразным изменением направления [3], [4], [5-8], [9-12], так и от моделей непрерывных случайных перемещений, с изначально детерминированной по величине, скоростью [13].

Цель данной работы, продемонстрировать особенности моделей с ортогональными воздействиями (1), (2), уравнений для определения распределения координаты частицы. Будут построены и уравнения, решение которых аппроксимируют распределение положения частицы, с учетом физической трактовки коэффициентов исходного стохастического уравнения.

**Замечание.** Исследование стохастических уравнений с ортогональными случайными влияниями, результат развития теории первых интегралов стохастических систем [14],[15]. Более высокая размерность пространства не является ограничением при построении таких моделей [15].

## 1. Особенности решения уравнения Ито с ортогональными воздействиями.

Рассмотрим вариант модели (1) изменения скорости в трехмерном пространстве:

$$dv(t) = -a(t)v(t)dt + [v(t) \times H]dt + \frac{b(t)}{|v(t)|}\left[v(t) \times dw(t)\right], \quad (3)$$

Предполагаем, что:

I. $a(t) > 0, \forall t \geq 0, a(t), b(t)$ - непрерывные и ограниченные $\forall t \geq 0$ скалярные функции; $v(t)$=$v(t; v(0))$, $v(0) \in R^3, v(0) \neq 0$;

II. $w(t) \in R^3$ – винеровский процесс с независимыми компонентами.

$[\cdot \times \cdot]$ – обозначение векторного произведения.

При ограничениях I-II, решение (3) существует и единственно [16].

**Замечание.** В общем случае, не обязательно требование независимости компонент $w(t)$[15].

Исследуем свойства $|\mathrm{v}(t)|^2$.

Воспользовавшись (3) и формулой Ито, находим:

$$d|\mathrm{v}(t)|^2 = [-2a(t)|\mathrm{v}(t)|^2 - 2b^2(t)]dt. \tag{4}$$

Т.е., модуль скорости - неслучайная величина.

Решение (4) имеет вид:

$$|\mathrm{v}(t)|^2 = \\ = 2\int_0^t d\theta b^2(\theta)\exp\left\{-2\int_\theta^t a(\tau)d\tau\right\} + \exp\left\{-2\int_0^t a(\tau)d\tau\right\}|\mathrm{v}(0)|^2 \tag{5}$$

Если $\dfrac{b^2(t)}{a(t)} = \dfrac{b^2(0)}{a(0)} = \dfrac{b^2}{a} = const, \forall t \geq 0$, то, после интегрирования первого слагаемого в (5) по частям, находим:

$$|\mathrm{v}(t)|^2 = \frac{b^2}{a}\left(1 - \exp\left\{-2\int_0^t a(\tau)d\tau\right\}\right) + \left\{-2\int_0^t a(\tau)d\tau\right\}|\mathrm{v}(0)|^2 \tag{6}$$

Из (6) следует, что $\forall \mathrm{v}(0) \neq 0$, $\lim\limits_{t\to\infty}|\mathrm{v}(t)|^2 = \dfrac{b^2}{a}$. .

Когда $a(t) = a,\ b(t) = b, \forall t \geq 0$ – постоянные, (6) приобретает более простой вид:

$$|\mathrm{v}(t)|^2 = \frac{b^2}{a}(1 - \exp\{-2at\}) + \exp\{-2at\}\,|\,\mathrm{v}(0)\,|^2 \tag{7}$$

Как видим, уравнение (3) моделирует релаксацию модуля скорости (уравнения (6), (7)), к постоянной величине. Т.е., сфера с радиусом $b(a)^{-0,5}$ - притягивающее многообразие [14], [15].

Перейдем к построению характеристической функции для процесса $x(t)$, зависящего от начальных значений $\{x(0); \mathrm{v}(0)\}$.

## 2. Уравнение для характеристической функции.

Для упрощения выкладок, в уравнении (3) положим $H = 0$, и рассмотри уравнение для составного процесса $\{x(t); \mathrm{v}(t)\}$ (в предположении независимости $x(0),\ \mathrm{v}(0)$):

$$\begin{cases} d\mathrm{v}(t) = -a(t)\mathrm{v}(t)dt + \dfrac{b(t)}{|\mathrm{v}(t)|}\left[\mathrm{v}(t) \times dw(t)\right], \\ x(t) = x(0) + \displaystyle\int_0^t \mathrm{v}(\tau)\, d\tau; \quad x(t), \mathrm{v}(t) \in \mathbb{R}^3. \end{cases} \tag{8}$$

Ограничения I-II уравнения (3), обеспечивают существование и единственность решение (8).

При построении уравнения для характеристической функции

$$\Psi_t = \mathrm{M}\left[\exp\{i(\lambda, x(t))\} / \mathrm{v}(0), x(0)\right],\ \Psi_0 = \exp\{i(\lambda, x(0))\}, \tag{9}$$

потребуется ряд вспомогательных утверждений, связанных с уравнениями для составного процесса $\{(\lambda, x(t)); (\lambda, \mathrm{v}(t))\}$.

Как следует из (8), с учетом свойства смешанного скалярно-векторного произведения,

$$\begin{cases} d(\lambda, \mathrm{v}(t)) = -a(t)(\lambda; \mathrm{v}(t))dt + \bar{b}(t)(dw(t), [\lambda \times \mathrm{v}(t)]), \\ d(\lambda, x(t)) = (\lambda, \mathrm{v}(t))dt \end{cases} \tag{10}$$

где $\bar{b}(t) = b(t) / |\mathrm{v}(t)|$.

**Лемма 1.**

$$d(\lambda, \mathrm{v}(t))^2 = -[2a(t) + \bar{b}^2(t)](\lambda, \mathrm{v}(t))^2 dt + \bar{b}^2(t)|\lambda|^2 |\mathrm{v}(t)|^2 dt + \\ + 2(\lambda, \mathrm{v}(t))\bar{b}(t)(dw(t), [\lambda \times \mathrm{v}(t)]) \tag{11}$$

**Доказательство.** Воспользовавшись (10) и формулой Ито, находим:

$$d(\lambda, \mathrm{v}(t))^2 = -2a(t)(\lambda, \mathrm{v}(t))^2 dt + \bar{b}^2(t)[(\lambda_2 \mathrm{v}_3(t) - \lambda_3 \mathrm{v}_2(t))^2 + \\ +(\lambda_1 \mathrm{v}_3(t) - \lambda_3 \mathrm{v}_1(t))^2 + (\lambda_1 \mathrm{v}_2(t) - \lambda_2 \mathrm{v}_1(t))^2]dt + \\ + 2\bar{b}(t)(\lambda, \mathrm{v}(t)_t)(dw(t), [\lambda \times \mathrm{v}(t)]) \tag{12}$$

Раскроем и преобразуем в (12) выражение в квадратных скобках при $dt$:

$$[(\lambda_2 \mathrm{v}_3(t) - \lambda_3 \mathrm{v}_2(t))^2 + (\lambda_1 \mathrm{v}_3(t) - \lambda_3 \mathrm{v}_1(t))^2 + (\lambda_1 \mathrm{v}_2(t) - \lambda_2 \mathrm{v}_1(t))^2] =$$
$$= \lambda_2^2 \mathrm{v}_3^2(t) - 2\lambda_2\lambda_3 \mathrm{v}_2(t)\mathrm{v}_3(t) + \lambda_3^2 \mathrm{v}_2^2(t) + \lambda_1^2 \mathrm{v}_3^2(t) -$$
$$-2\lambda_1\lambda_3 \mathrm{v}_1(t)\mathrm{v}_3(t) + \lambda_3^2 \mathrm{v}(t)_1^2 + \lambda_1^2 \mathrm{v}_2^2(t) - 2\lambda_1\lambda_2 \mathrm{v}_1(t)\mathrm{v}_2(t) + \lambda_2^2 \mathrm{v}_1^2(t) =$$
$$= \lambda_1^2 |\mathrm{v}(t)|^2 + \lambda_2^2 |\mathrm{v}(t)_t|^2 + \lambda_3^2 |\mathrm{v}(t)|^2 -$$
$$-(\lambda_1^2 \mathrm{v}_1^2(t) + \lambda_2^2 \mathrm{v}_2^2(t) + \lambda_3^2 \mathrm{v}_3^2(t) + 2\lambda_2\lambda_3 \mathrm{v}_2(t)\mathrm{v}_3(t) +$$
$$+ 2\lambda_1\lambda_3 \mathrm{v}_1(t)\mathrm{v}_3(t) + 2\lambda_1\lambda_2 \mathrm{v}_1(t)\mathrm{v}_2(t)) =$$
$$= |\lambda|^2 |\mathrm{v}(t)|^2 - (\lambda_1 \mathrm{v}_1(t) + \lambda_2 \mathrm{v}(t)_2 +$$
$$+\lambda_3 \mathrm{v}_3(t))(\lambda_1 \mathrm{v}_1(t) + \lambda_2 \mathrm{v}_2(t) + \lambda_3 \mathrm{v}_3(t)) = |\lambda|^2 |\mathrm{v}(t)|^2 - (\lambda, \mathrm{v}(t))^2$$

Учитывая это преобразование, возвратимся к уравнению (12):

$$d(\lambda, \mathrm{v}(t))^2 = -2a(t)(\lambda, \mathrm{v}(t))^2 dt + \\ +\bar{b}^2(t)[|\lambda|^2 |\mathrm{v}(t)|^2 - (\lambda, \mathrm{v}(t))^2]dt + 2\bar{b}(t)(\lambda, \mathrm{v}(t))(dw(t), [\lambda \times \mathrm{v}(t)]) \tag{13}$$

После перегруппировки слагаемых в (13), приходим к уравнению (11) .

**Лемма 2.**

$$\frac{\partial \Psi_t}{\partial t} + a(t)\frac{\partial^2 \Psi_t}{\partial t^2} = \mathrm{M}\left[(\lambda, \mathrm{v}(t))^2 \exp\{i(\lambda, x(t))\} / \mathrm{v}(0); x(0)\right]. \tag{14}$$

**Доказательство.** Воспользовавшись определением (9) и последовательно дважды формулой Ито, находим:

$$\frac{\partial \Psi_t}{\partial t} = \mathrm{M}\left[i(\lambda, \mathrm{v}(t))\exp\{i(\lambda, x(t))\} / \mathrm{v}(0); x(0)\right], \tag{15}$$

$$\frac{\partial^2 \Psi_t}{\partial t^2} dt = d\frac{\partial \Psi_t}{\partial t} = \mathrm{M}\left[id[(\lambda, \mathrm{v}(t))\exp\{i(\lambda, x(t))\}] / \mathrm{v}(0); x(0)\right] =$$
$$= dt\mathrm{M}\left[[-(\lambda, \mathrm{v}(t))^2 dt + i(\lambda, d\mathrm{v}(t))]\exp\{i(\lambda, x(t))\} / \mathrm{v}(0); x(0)\right] =$$
$$= dt\mathrm{M}\left[-([\lambda, \mathrm{v}(t))^2 dt - ia(t)(\lambda, \mathrm{v}(t))dt]\times \right. \tag{16}$$
$$\left. \times\exp\{i(\lambda, x(t))\} / \mathrm{v}(0); x(0)\right] =$$
$$= -dt\mathrm{M}\left[(\lambda, \mathrm{v}(t))^2 \exp\{i(\lambda, x(t))\} / \mathrm{v}(0); x(0)\right] - a(t)\frac{\partial \Psi_t}{\partial t} dt,$$

Тогда, с учетом (15), уравнения для второй производной (16) приходим к равенству (14).

**Лемма 3.**

$$
\begin{aligned}
&-\mathrm{M}[(\lambda,\mathrm{v}(t))^2\exp\{i(\lambda,x(t))\}/\mathrm{v}(0);\,x(0)]=\\
&=\mathrm{M}[\exp\{i(\lambda,x(t))\}/\mathrm{v}(0);\,x(0)]\times\\
&\times\left[-(\lambda,\mathrm{v}(0))^2-\int_0^t \bar{b}^2(\tau)|\lambda|^2|\mathrm{v}(\tau)|^2\exp\{-\int_\tau^t[2a+\bar{b}^2(\theta)]d\theta\}d\tau\right]
\end{aligned} \quad (17).
$$

**Доказательство**. Представим решение уравнения (11) в таком виде:

$$
\begin{aligned}
&(\lambda,\mathrm{v}(T-t;\mathrm{v}(t)))^2=\exp\{-\int_t^T[2a(\theta)+\bar{b}^2(\theta)]d\theta\}(\lambda,\mathrm{v}(t))^2+\\
&+\int_t^T\exp\{-\int_\tau^T[2a(\theta)+\bar{b}^2(\theta)]d\theta\}b^2(\tau)|\lambda|^2|\mathrm{v}(\tau)|^2\,d\tau+\\
&+\int_t^T\exp\{-\int_\tau^T[2a(\theta)+\bar{b}^2(\theta)]d\theta\}2(\lambda,\mathrm{v}(\tau))\bar{b}(\tau)(dw(\tau),[\lambda\times\mathrm{v}(\tau)])
\end{aligned} \quad . \quad (18)
$$

Из условий существование и единственность решения (11) , следует что $\mathrm{v}(T-t;\mathrm{v}(t))=\mathrm{v}(T)$. Учтем это и перейдем от равенства (18) к такому:

$$
\begin{aligned}
&(\lambda,\mathrm{v}(t))^2=\exp\{-\int_0^t[2a(\theta)+\bar{b}^2(\theta)]d\theta\}\times\\
&\times\Bigg[(\lambda,\mathrm{v}(T)\,)^2\exp\{\int_0^T[2a(\theta)+\bar{b}^2(\theta)]d\theta\}-\\
&-\int_0^T\exp\{\int_0^\tau[2a(\theta)+\bar{b}^2(\theta)]d\theta\}b^2(\tau)|\lambda|^2|\mathrm{v}(\tau)|^2\,d\tau-\\
&-\int_t^T\exp\{\int_0^\tau[2a(\theta)+\bar{b}^2(\theta)]d\theta\}2(\lambda,\mathrm{v}(\tau))\bar{b}(\tau)(dw(\tau),[\lambda\times\mathrm{v}(\tau)])\Bigg]+\\
&+\int_0^t\bar{b}^2(\tau)|\lambda|^2|\mathrm{v}(\tau)|^2\exp\{-\int_\tau^t[2a(\theta)+\bar{b}^2(\theta)]d\theta\}d\tau
\end{aligned} \quad (19)
$$

Опираясь на (19), преобразуем $M[(\lambda,\mathrm{v}(t))^2\exp\{i(\lambda,x(t))\}/x(0),v(0)]$, с учетом детерминированности $|\mathrm{v}(\tau)|$ (см. (4)):

$$
\begin{aligned}
&|-\mathrm{M}[(\lambda,\mathrm{v}(t))^2\exp\{i(\lambda,x(t))\}/x(0);\mathrm{v}(0)]=-\exp\{\int_0^t[2a(\theta)+\bar{b}^2(\theta)]d\theta\}\times\\
&\times\mathrm{M}[\exp\{i(\lambda,x(t))\}\Bigg[(\lambda,\mathrm{v}(T))^2\exp\{+\int_0^T[2a(\theta)+\bar{b}^2(\theta)]d\theta\}-\\
&-\int_0^T\exp\{\int_0^\tau[2a(\theta)+\bar{b}^2(\theta)]d\theta\}b^2(\tau)|\lambda|^2|\mathrm{v}(\tau)|^2\,d\tau\Bigg]/\mathrm{v}(0);x(0))]-\\
&-\mathrm{M}\Big[\exp\{i(\lambda,x(t))\}\times\\
&\times\int_t^T\exp\{-\int_\tau^t[2a(\theta)+\bar{b}^2(\theta)]d\theta\}2(\lambda,\mathrm{v}(\tau))\bar{b}(\tau)(dw(\tau),[\lambda\times\mathrm{v}(\tau)])/\mathrm{v}(0);\,x(0))\Big]
\end{aligned} \quad (20)
$$

Выражение в (20) слева не зависит от $T$. Поэтому, для первого слагаемого справа в (20), $\forall t\in(T,0]$, $\forall T\geq 0$, должно выполняться равенство:

$$-\mathrm{M}[\exp\{i(\lambda,x(t))\}\Bigg[(\lambda,\mathrm{v}(T))^2\exp\{\int_0^T[2a(\theta)+\bar{b}^2(\theta)]d\theta\}-$$

$$-\int_0^T\exp\{\int_0^\tau[2a(\theta)+\bar{b}^2(\theta)]d\theta\}\bar{b}^2(\tau)|\lambda|^2|\mathrm{v}(\tau)|^2d\tau\Bigg]/\mathrm{v}(0);x(0)]= \qquad (21)$$

$$=-\mathrm{M}[\exp\{i(\lambda,x(t))\}/\mathrm{v}(0);\ x(0)](\lambda,\mathrm{v}(0))^2.$$

Далее, поскольку $w(\tau)$- процесс с независимыми приращениями, т.е., $\forall\delta>0$ - $[w(\tau+\delta)-w(\tau)]$ не зависят от значений $x(\tau),\mathrm{v}(\tau)$, и, поскольку, $\mathrm{M}[w(\delta+\tau)-w(\delta)]=0$, приходим к равенству для второго слагаемого в (20):

$$\mathrm{M}\Big[\exp\{i(\lambda,x(t))\}\times$$

$$\times\int_t^T\exp\{-\int_\tau^t[2a(\theta+\bar{b}^2(\theta)]d\theta\}2(\lambda,\mathrm{v}(\tau))\bar{b}(\tau)\times \qquad (22)$$

$$\times(dw(\tau),[\lambda\times\mathrm{v}(\tau)])/x(0);\mathrm{v}(0))\Big]\equiv 0$$

Проверить это равенство можно и непосредственно, взяв производную по $T$, воспользовавшись формулой Ито и уравнением (11).

С учетом (21) и (22), уравнения (20) переходит в уравнение (17).

**Теорема.** Характеристическая функция

$$\Psi_t=\mathrm{M}\Big[\exp\{i(\lambda,x(t))\}/\mathrm{v}(0);x(0)\Big],\ \forall\mathrm{v}(0)\neq 0,$$

для процесса $x(t)$, $\lambda\in\mathbb{R}^3$, подчиненного системе (8), является решением уравнения:

$$\frac{\partial^2\Psi_t}{\partial t^2}+a(t)\frac{\partial\Psi_t}{\partial t}=-|\lambda|^2\Psi_t f(t)-(\lambda,\mathrm{v}(0))^2\Psi_t \qquad (23.а)$$

$$\Psi_0=\exp\{i(\lambda,x(0))\},\frac{\partial\Psi_t}{\partial t}\Big|_{t=0}=(\lambda,\mathrm{v}(0))\Psi_0, \qquad (23.в)$$

где

$$f(t)=\int_0^t\bar{b}^2(\tau)|\mathrm{v}(\tau)|^2\exp\{-\int_\tau^t[2a(\theta)+\bar{b}^2(\theta)]d\theta\}d\tau]=$$

$$=\int_0^t\frac{b^2(\tau)}{|\mathrm{v}(\tau)|^2}|\mathrm{v}(\tau)|^2\exp\{-\int_\tau^t[2a(\theta)+\bar{b}^2(\theta)]d\theta\}d\tau]=$$

$$=\int_0^t b^2(\tau)\exp\{-\int_\tau^t 2a(\theta)[1+(2a(\theta))^{-1}|\mathrm{v}(\theta)|^{-2}b^2(\theta)]d\theta\}d\tau]\geq 0$$

**Доказательство Теоремы** – есть следствие замены выражения справа в (14) (Лемма 2) на выражение (17) (Лемма 3), определения (9) и равенства (15):

$$\frac{\partial\Psi_t}{\partial t}\Big|_{t=0}=$$
$$=\mathrm{M}\Big[i(\lambda,\mathrm{v}(t))\exp\{i(\lambda;x(t))\}/\mathrm{v}(0);\ x(0)\Big]\Big|_{t=0}=i(\lambda,\mathrm{v}(0))\exp\{i(\lambda;x(0))\} \qquad (24)$$

### 3. **Уравнение для плотности вероятности координаты.**

Перейдем, с учетом связи между характеристическими функциями и отвечающих им распределениям, от уравнения (23) к уравнению в исходных переменных по $x$:

$$\frac{\partial^2 \rho(t;x/x(0);\mathrm{v}(0))}{\partial t^2}+a(t)\frac{\partial \rho(t;x/x(0);\mathrm{v}(0))}{\partial t}=f(t)\nabla_x^2\rho(t;x/x(0);\mathrm{v}(0)).+$$

$$+\sum_{l=1}^{3}\sum_{j=1}^{3}\mathrm{v}_l(0)\mathrm{v}_j(0)\frac{\partial^2 \rho(t;x/x(0);\mathrm{v}(0))}{\partial x_l \partial x_j} \qquad (25)$$

Уравнение для плотности распределения

$$\rho(t;x/\mathrm{v}(0))=\int_{-\infty}^{\infty}\rho(t;x/\acute{o};\mathrm{v}(0))\rho(\acute{o})d\acute{o}_1 d\acute{o}_2 d\acute{o}_3,$$

соответственно будет таким:

$$\frac{\partial^2 \rho(t;x/\mathrm{v}(0))}{\partial t^2}+a(t)\frac{\partial\ \rho(t;x/;\mathrm{v}(0))}{\partial t}=f(t)\nabla_x^2\rho(t;x/\mathrm{v}(0)).+$$

$$+\sum_{l=1}^{3}\sum_{j=1}^{3}\mathrm{v}_l(0)\mathrm{v}_j(0)\frac{\partial^2 \rho(t;x/\mathrm{v}(0))}{\partial x_l \partial x_j}.$$

Рассмотрим случай когда $a,\ b,\ \left|\mathrm{v}(t)\right|^2=\frac{b^2}{a}=\left|\mathrm{v}\right|^2$ - постоянные. Тогда

$$f(t)=\int_0^t b^2 \mathrm{e\,xp}\{-3a(t-\tau)\}d\tau]=\frac{b^2}{3a}[1-\mathrm{e\,xp}\{-3at\}]=3^{-1}|\,\mathrm{v}\,|^2[1-\mathrm{e\,xp}\{-3at\}] \qquad (26)$$

поскольку, как следует из (7), $a^{-1}b^2=|\,\mathrm{v}\,|^2$ - равновесное значение.

В этом случае, (25) приобретает вид:

$$\frac{\partial^2 \rho(t;x/\mathrm{v}(0))}{\partial t^2}+a\frac{\partial\ \rho(t;x/\mathrm{v}(0))}{\partial t}=$$

$$=3^{-1}|\,\mathrm{v}\,|^2[1-\mathrm{e\,xp}\{-3at\}]\nabla_x^2\rho(t;x/\mathrm{v}(0))+ \qquad (27)$$

$$+\sum_{l=1}^{3}\sum_{j=1}^{3}\mathrm{v}_l(0)\mathrm{v}_j(0)\frac{\partial^2 \rho(t;x/\mathrm{v}(0))}{\partial x_l \partial x_j},$$

Решение (27), в соответствии с (23.в), ищем при начальных условиях:

$$\rho(t;x/\mathrm{v}(0))|_{t=0}=\rho(x);\quad \frac{\partial \rho(t;x/\mathrm{v}(0))}{\partial t}|_{t=0}=\left(\mathrm{v}(0),\vec{\nabla}_x\rho(x)\right) \qquad (28)$$

## 4. Уравнения для аппроксимирующих решений.

Рассмотрим варианты (27) которые имеют физическую трактовку, и построим уравнения, решения которых можно рассматривать как аппроксимирующие решения уравнения (27).

### 4.1. Диффузионное приближение

Положим, что в уравнении (8)

$$a=\tilde{a}\varepsilon^{-1},\ \ b=\tilde{b}\varepsilon^{-1}, \qquad (29)$$

где и $\tilde{a},\ \tilde{b}$ ограниченные величины.

В данном случае,

$$\tilde{b}^2\varepsilon^{-2}/\varepsilon^{-1}\tilde{a}=\varepsilon^{-1}\tilde{b}^2/\tilde{a}=\varepsilon^{-1}|\,\tilde{\mathrm{v}}\,|^2.$$

При условии (29), уравнении (27) приобретает вид:

$$\varepsilon\frac{\partial^2 \rho_\varepsilon(t;x/\mathrm{v}(0))}{\partial t^2}+\tilde{a}\frac{\partial \rho_\varepsilon(t;x/\mathrm{v}(0))}{\partial t}=$$

$$=3^{-1}|\,\tilde{\mathrm{v}}\,|^2[1-\mathrm{e\,xp}\{-3\varepsilon^{-1}\tilde{a}t\}]\nabla_x^2\rho_\varepsilon(t;x/\mathrm{v}(0))+ \qquad (30)$$

$$+\sum_{l=1}^{3}\sum_{j=1}^{3}\tilde{\mathrm{v}}_l(0)\tilde{\mathrm{v}}_j(0)\frac{\partial^2 \rho_\varepsilon(t;x/\mathrm{v}(0))}{\partial x_l \partial x_j},\ \tilde{\mathrm{v}}_l(0)\frac{1}{\sqrt{\varepsilon}}=\mathrm{v}(0)$$

**Замечание.** Выбор коэффициентов в таком виде связан с тем, что как показано, например, в [17], [18], при устремлении $\varepsilon$ к нулю, приходим к правильным уравнения для диффузионного приближения положения броуновской частицы.

Другим моментом является то, что при переходе к безразмерным величинам в уравнении Ланжевена, для комнатных температур, коэффициенты приобрели вид (29) с $\varepsilon \sim 10^{-12}$[15]. Последующий переход от этого уравнения к диффузионному приближению по координате, привел к уравнению для плотности частиц, совпадающим со вторым законом Фика (уравнением диффузии для однородных и неоднородных сред) [15]. Это дает нам основание, исследовать уравнение (3), с коэффициентами (29) и для модели с ортогональными воздействиями.

Если мы заменим в (27), (30) $\mathrm{v}(0) \to -\mathrm{v}(0)$, то вид этих уравнений останется прежним. Т.е., при таком изменении направления $\mathrm{v}(0)$

$$\rho_\varepsilon(t; x/\mathrm{v}(0)) = \rho_\varepsilon(t; x/-\mathrm{v}(0)) \ . \qquad (31)$$

В соответствии с выводами работ [14], [15], направления случайного вектора $\tilde{\mathrm{v}}(\infty)$ равномерно распределены на сфере.

Выбирая в качестве начального условия это требование равномерного распределение на сфере, приходим к равенствам:

$$|\mathrm{v}(t)|^2 = |\mathrm{v}(0)|^2 =$$

$$= \mathrm{M}[|\mathrm{v}_1(t)|^2/\mathrm{v}(0)] + \mathrm{M}[|\mathrm{v}_2(t)|^2/\mathrm{v}(0)] + \mathrm{M}[|\mathrm{v}_3(t)|^2/\mathrm{v}(0)], \ \forall t \ge 0.$$

Следовательно, в силу равноправности всех проекций, с учетом того, что $|\mathrm{v}(t)|^2 = |\mathrm{v}|^2 = const$, должно выполняться равенство:

$$\mathrm{M}[|\mathrm{v}_1(t)|^2/\mathrm{v}(0)] = \mathrm{M}[|\mathrm{v}_2(t)|^2/\mathrm{v}(0)] = \mathrm{M}[|\mathrm{v}_3(t)|^2/\mathrm{v}(0) = |\mathrm{v}|^2\, 3^{-1}$$

Поэтому, в сферической системе координат, необходимо выбрать представление для проекции в исходном распределении в таком виде:

$$\tilde{\mathrm{v}}_z(0) = 3^{-0,5}|\tilde{\mathrm{v}}|\cos\theta, \tilde{\mathrm{v}}_x(0) = 3^{-0,5}|\tilde{\mathrm{v}}|\sin\theta\cos\varphi, \tilde{\mathrm{v}}_y(0) = 3^{-0,5}|\tilde{\mathrm{v}}|\sin\theta\sin\varphi.$$

Рассмотрим выражение:

$$\frac{1}{(2\pi)^2}\int_{-\pi}^{\pi} d\theta \int_{-\pi}^{\pi} \tilde{\mathrm{v}}_l(0;\varphi;\theta)\tilde{\mathrm{v}}_j(0;\varphi;\theta)\rho_\varepsilon(t; x/\tilde{\mathrm{v}}(0))d\varphi =$$

$$= \frac{|\tilde{\mathrm{v}}|^2}{(2\pi)^2}\int_{-\pi}^{\pi} d\theta \int_{-\pi}^{\pi} d\varphi((\sin\theta)^2(\sin\varphi\cos\varphi)\times \qquad (32)$$

$$\times\rho_\varepsilon(t; x/|\tilde{\mathrm{v}}|\sin\theta\sin\varphi; |\tilde{\mathrm{v}}|\sin\theta\cos\varphi; |\tilde{\mathrm{v}}|\cos\theta) = 0, \forall j \ne l.$$

в силу свойства (31).

Рассмотрим, теперь, выражение:

$$M[\tilde{\mathrm{v}}_l^2(0)/x;t] =$$

$$= \frac{1}{(2\pi)^2\rho_\varepsilon(t;x)]}\int_{-\pi}^{\pi} d\theta \int_{-\pi}^{\pi} \tilde{\mathrm{v}}_l^2(0;\varphi;\theta)\rho_\varepsilon(t; x/\tilde{\mathrm{v}}(0))d\varphi \ne 0, \ \ l = 1,2,3.$$

Из условия сферической симметрии

$$[(\tilde{\mathrm{v}}_1(0))^2/x;t] = [(\tilde{\mathrm{v}}_2(0))^2/x;t] = [(\tilde{\mathrm{v}}_3(0))^2/x;t]$$

и требования

$$[(\tilde{\mathrm{v}}_1(0))^2/x;t] + [(\tilde{\mathrm{v}}_2(0))^2/x;t] + [(\tilde{\mathrm{v}}_3(0))^2/x;t] = |\tilde{\mathrm{v}}(t)|^2 = |\tilde{\mathrm{v}}|^2,$$

следует, что

$$[(\tilde{\mathrm{v}}_l(0))^2/x;t] = \frac{1}{3}|\tilde{\mathrm{v}}|^2, \forall l - 1,2,3, \forall x; \forall t \ge 0 \ . \qquad (33)$$

С четом равенств (32), (33), после усреднения по начальному распределению $\tilde{v}(0)$, приходим к уравнению:

$$\varepsilon\frac{\partial^2\rho_\varepsilon(t;x)}{\partial t^2}+\tilde{a}\frac{\partial\rho_\varepsilon(t;x)}{\partial t}= \\ =3^{-1}|\tilde{v}|^2[1-\exp\{-3\varepsilon^{-1}\tilde{a}t\}]\nabla_x^2\rho_\varepsilon(t;x)+3^{-1}|\tilde{v}|^2\nabla_x^2\rho_\varepsilon(t;x) \tag{34}$$

Как было доказано в [15] (Гл. 1,3.2, Теорема 3.3), когда $\varepsilon\to 0$ процесс $x_\varepsilon(t)$, для произвольных начальных распределений скорости частиц, при ограничениях на четвертые моменты: $M|v(0)|^4\le\varepsilon^{-2}K, K=const$ - слабо сходится к решению стохастического уравнения:

$$dx(t)=\frac{\tilde{b}}{\tilde{a}}\sqrt{\frac{2}{3}}dw(t)$$

Уравнение для плотности распределение для этого процесса имеет вид:

$$\frac{\partial\ \rho(t;x)}{\partial t}=\frac{1}{3}|\tilde{v}|^2\ \tilde{a}^{-1}\nabla_x^2\rho(t;x)$$

- и совпадает с (34), если отбросить в (34) слагаемые с $\varepsilon$ при устремлении $\varepsilon$ к 0.

Проведенные преобразования и сопоставления, служат нам основанием, при условиях (29), использовать, в качестве аппроксимирующего для решения (27), решение уравнения:

$$a_0\frac{\partial\ \bar{\rho}(t;x/\tilde{v}(0))}{\partial t}= \\ =3^{-1}|\tilde{v}|^2\ \nabla_x^2\bar{\rho}(t;x/\tilde{v}(0))+\sum_{l=1}^{3}\sum_{j=1}^{3}\tilde{v}_l(0)\tilde{v}_j(0)\frac{\partial^2\bar{\rho}(t;x/\tilde{v}(0))}{\partial x_l\partial x_j}. \tag{35}$$

Выберем оси координат таким образом, чтобы направление оси $x_1$ совпадало с направлением $\pm v(0)$. В этом случае, приходим к такому уравнению для плотности распределения $\bar{\rho}(t;x/\tilde{v}(0))$:

$$\frac{\partial\bar{\rho}(t;x/\tilde{v}_1(0))}{\partial t}=\frac{|\tilde{v}|^2}{3a}\nabla_x^2\bar{\rho}(t;x/\tilde{v}_1(0))+\frac{|\tilde{v}|^2}{3a}\frac{\partial^2\bar{\rho}(t;x/\tilde{v}_1(0))}{\partial x_1^2}, \tag{36}$$

т.к. $v_2(0)=0, v_3(0)=0\Rightarrow|v|^2=|v_1|^2$.

После выполнения замены:

$$u_1=\sqrt{\frac{3a}{2|\tilde{v}|^2}}x_1;\ u_2=\sqrt{\frac{3a}{|\tilde{v}|^2}}x_2;\ u_3=\sqrt{\frac{3a}{|\tilde{v}|^2}}x_3,$$

(36) переходит в обыкновенное уравнение диффузии:

$$\frac{\partial\bar{\bar{\rho}}(t;u/\tilde{v}_1(0))}{\partial t}=\nabla_u^2\bar{\bar{\rho}}(t;u/\tilde{v}_1(0))\ , \tag{37}$$

где $\bar{\bar{\rho}}(t;u/\tilde{v}_1(0))=\bar{\rho}(t;u_1(\frac{3a}{2|\tilde{v}|^2})^{-0,5};u_2(\frac{3a}{|\tilde{v}|^2})^{-0,5};u_3(\frac{3a}{|\tilde{v}|^2})^{-0,5}/\tilde{v}_1(0))$

Решение (37) имеет вид:

$$\bar{\bar{\rho}}(t;u/\tilde{v}(0))= \\ =\int_{-\infty}^{\infty}\frac{1}{2\sqrt{\pi t}}\exp\left(-\frac{[(u_1-y_1)^2+(u_2-y_2)^2+(u_3-y_3)^2}{4t}\right)\bar{\bar{\rho}}(y/\tilde{v}_1(0))dy_1dy_2dy_3.$$

Возвратившись к исходным переменным, находим решение уравнения (36):

$$\bar{\rho}(t;x/\tilde{\mathrm{v}}_1(0)) =$$

$$= \int_{-\infty}^{\infty}\left(\frac{3a}{|\tilde{\mathrm{v}}|^2}\right)^{\frac{3}{2}}\frac{1}{2\sqrt{2\pi t}}\exp\left(-\frac{3a}{2|\tilde{\mathrm{v}}|^2}\frac{(x_1-z_1)^2}{4t}-\frac{3a}{|\tilde{\mathrm{v}}|^2}\frac{[(x_2-z_2)^2+(x_3-z_3)^2]}{4t}\right)\times$$

$$\times\rho(z/\tilde{\mathrm{v}}_1(0))dz_1dz_2dz_3$$

Как видим, наблюдается ассиметрия расплывания плотности: более медленное по направлению $\pm\tilde{\mathrm{v}}(0)$.

### 4.2. Случай малого стоксового трения и малых возмущений. Переход к волновому уравнению

Пусть в уравнении (8),

$$a=\varepsilon a_0 \text{ а } b^2/\varepsilon a_0 = |\mathrm{v}|^2 = const$$

Для согласований этих требований необходимо положить, что

$$b^2/a_0 = \varepsilon a_0|\mathrm{v}|^2, \text{ т.е., } b=b_0\sqrt{\varepsilon},\ b_0=|\mathrm{v}|\sqrt{\varepsilon a_0}.$$

Этому случаю будет соответствовать такой вид уравнения (27):

$$\frac{\partial^2\rho_\varepsilon(t;x/\mathrm{v}(0))}{\partial t^2}+\varepsilon a\frac{\partial\ \rho_\varepsilon(t;x/\mathrm{v}(0))}{\partial t}=$$

$$=3^{-1}|\mathrm{v}|^2[1-\exp\{-3a_0\varepsilon t\}]\nabla_x^2\rho_\varepsilon(t;x/\mathrm{v}(0))+ \qquad (38)$$

$$+\sum_{l=1}^{3}\sum_{j=1}^{3}\mathrm{v}_l(0)\mathrm{v}_j(0)\frac{\partial^2\rho_\varepsilon(t;x/\mathrm{v}(0))}{\partial x_l\partial x_j},$$

На временных интервалах, когда $t<\dfrac{1}{\varepsilon^\alpha a_0}, 0<\alpha<1$, $([1-\exp\{-3\varepsilon a_0 t\}]=O(\varepsilon^{1-\alpha}))$, в качестве аппроксимирующего, выбираем решение уравнения:

$$\frac{\partial^2\rho(t;x/\mathrm{v}(0))}{\partial t^2}=\sum_{l=1}^{3}\sum_{j=1}^{3}\mathrm{v}_l(0)\mathrm{v}_j(0)\frac{\partial^2\rho(t;x/\mathrm{v}(0))}{\partial x_l\partial x_j}, \qquad (39)$$

Уравнение (38) приобретает более простой вид, если одна из осей выбранной системы координат, например $x_1$, совпадает с + или - направлением скорости $\mathrm{v}(0)$.

В этом случае, уравнение (39) заменяется таким:

$$\frac{\partial^2\tilde{\rho}(\tau;z/\mathrm{v}_1(0))}{\partial t^2}=|\mathrm{v}|^2\frac{\partial^2\tilde{\rho}(\tau;z/\mathrm{v}_1(0))}{\partial x_1^2}, \qquad (40)$$

т.к., $\mathrm{v}_2(0)=0, \mathrm{v}_3(0)=0 \Rightarrow |\mathrm{v}|^2=|\mathrm{v}_1|^2$.

Это волновое уравнение. Его решение - произвольная функция

$$\rho(t;x_1;x_2;x_3/\mathrm{v}_1(0))=q_1([x_1-t|\mathrm{v}_1(0)|];x_2;x_3)+q_2([x_1+t|\mathrm{v}_1(0)|];x_2;x_3),$$
$$\rho(0;x_1;x_2;x_3/\mathrm{v}_1(0))=\rho(x_1;x_2;x_3). \qquad (41)$$

Из (28) следует, что при $t=0$,

$$\frac{\partial\rho(t;x/\mathrm{v}_1(0))}{\partial t}\Big|_{t=0}=\mathrm{v}_1(0)\frac{\partial\rho(x)}{\partial x_1},$$

Подставив представление для (41) в это равенство, находим:

$$\frac{\partial}{\partial t}q_1([x_1-t|\mathrm{v}_1(0)|];x_2;x_3)+\frac{\partial}{\partial t}q_2([x_1+t|\mathrm{v}_1(0)|];x_2;x_3)=$$

$$=-|\mathrm{v}_1|\frac{\partial}{\partial x_1}q_1(x_1;x_2;x_3)+|\mathrm{v}_1|\frac{\partial}{\partial x_1}q_2(x_1;x_2;x_3)=$$

$$=\mathrm{v}_1\frac{\partial}{\partial x_1}q_1(x_1;x_2;x_3)+\mathrm{v}_1\frac{\partial}{\partial x_1}q_2(x_1;x_2;x_3)$$

Откуда следует, что если $v_1(0)$ направлено в положительном направлении оси $x_1$, то $q_1(x_1;x_2;x_3)=0$, и $q_2(x_1;x_2;x_3)=\rho(x_1;x_2;x_3)$. Если же направление $v_1(0)$ противоположно направлению оси $x_1$, то $q_2(x_1;x_2;x_3)=0$, и $q_1(x_1;x_2;x_3)=\rho(x_1;x_2;x_3)$.

## Заключение.

Изменения реальных процессов, реальных сред, в том числе и стохастических, всегда связаны с ограничениями на их вариации [19]. Перемещения не происходят с бесконечно большой скоростью. В то же время, классические уравнения диффузии основаны явно на предположении, что скорости броуновских частиц могут принимать сколь угодно большие значения. В данной работе, на основе уравнений Ланжевена, при ортогональных к скорости частицы случайных воздействиях, продемонстрировано, что для таких систем могут существовать притягивающие многообразия, не смотря на то, что винеровский процесс – процесс, допускающий сколь угодно большие значения. В отличие от предыдущих работ [14], [15] в статье построено уравнение для определения вероятности распределения диффундирующей частицы в пространстве координат, включающую зависимость от вектора начальной скорости. Показано, что малые воздействия, при определенном согласовании коэффициентов в исходном стохастическом уравнении, приводят к описанию перемещения частиц на основе волнового уравнениям с постоянной скоростью.

Рассмотренные уравнения не исчерпывают класс моделей с ортогональными возмущениями. Более полный класс рассмотрен в [14], [15], в том числе, и для $n$-мерных процессов [15], в связи с разработкой теории первых интегралов для стохастических систем.

## Литература